\documentclass[]{ensmath}

\usepackage{epsfig}

\def\limi{\lim_{i\to\infty}}
\def\limn{\lim_{n\to\infty}}

\def\om{\omega}
\def\e{\epsilon}

\def\cGG{\cG=\{G_n\}^\infty_{n=1}}
\def\cHH{\cH=\{H_n\}^\infty_{n=1}}

\def\deg{\mbox{deg}\,}
\def\ke{\mbox{Ker}\,}
\def\im{\mbox{Im}\,}
\def\rank{\mbox{Rank}\,}

\def\dim{{\rm dim}}

\def\Fo{F$\mbox{\o}$lner}

\def\proof{\smallskip\noindent{\it Proof.} }
\def\Fp{{\mathbb F}_p}

\def\bZ{{\mathbb Z}}

\def\bQ{{\mathbb Q}}
\def\deg{\mbox{deg}\,}

\def\cF{\mbox{$\cal F$}}

\def\cG{\mbox{$\cal G$}}
\def\cH{\mbox{$\cal H$}}
\def\cJ{\mbox{$\cal J$}}
\def\cT{\mbox{$\cal T$}}

\def\to{\rightarrow}

\def\questionname{Question}
\newtheorem{question}[theorem]{\questionname}
\EnsMath pages 1-2
\title[THE COMBINATORIAL COST]{THE COMBINATORIAL COST}

\author[G. ELEK]{G\'abor \sn{Elek}\thanks{The 
author is supported by OTKA Grants T 049841 and T 037846 }}

\begin{document}

\maketitle

\begin{abstract} We study the combinatorial analogues of the classical
invariants of measurable equivalence relations. We introduce the notion
of cost and $\beta$-invariants (the analogue of the first $L^2$-Betti number
introduced by Gaboriau \cite{Gab2})
for sequences of finite graphs with uniformly bounded vertex degrees and
examine the relation of these invariants and the rank gradient resp.
mod $p$ homology gradient invariants
 introduced by Lackenby \cite{Lac1} \cite{Lac2} for residually finite groups.

\end{abstract}

\section{Introduction}
\subsection{Graph sequences}
Let $\cG=\{G_n\}^\infty_{n=1}$ be a sequence of finite simple
graphs satisfying the following conditions:
\begin{itemize}
\item $\sup_{1\leq n <\infty} \max_{x\in V(G_n)} \deg(x)<\infty $.
That is the graphs have uniformly bounded vertex degrees.
\item $|V(G_n)|\to \infty$ as $n\to\infty$.
\end{itemize}
In the sequel we refer to such systems as graph sequences. Now let
$\cH=\{H_n\}^\infty_{n=1}$ be another graph sequence such that
$V(H_n)=V(G_n)$ for any $n\geq 1$.
Then $\cH\prec \cG$ is there exists an integer $L>0$ such that
for any $n\geq 1$ and $x,y\in V(H_n)$, $d_{G_n}(x,y)\leq L d_{H_n}(x,y)$,
where $d_{G_n}$ resp. $d_{H_n}$ denote the shortest path metrics on
$G_n$ resp. on $H_n$. That is if $x$ and $y$ are adjacent in the graph
$H_n$ then there exists a path between $x$ and $y$ in $G_n$ of length
at most $L$. We say that $\cG$ and $\cH$ are {\bf equivalent}, $\cG\simeq \cH$
if $\cH\prec\cG$ and $\cG\prec\cH$. The {\bf edge number} of $\cG$ is
defined as
$$e(\cG):=\liminf_{n\to\infty} \frac{|E(G_n)|}{|V(G_n)|}$$
and the {\bf cost} of $\cG$ is given as
$$c(\cG):=\inf_{\cH\simeq \cG} e(\cH)\,.$$
Clearly, $c(\cG)\geq 1$ for any graph sequence $\cG$. Originally, the
cost was defined for measurable equivalence relations by Levitt \cite{Lev}.
In our paper we view graph sequences as the analogues of $L$-graphings
of measurable equivalence relations (see \cite{KM}). 

\noindent
Recall that a graph sequence $\cG=\{G_n\}^\infty_{n=1}$ is a large girth
sequence if for any $k\geq 1$, there exists $n_k$ such that if $n\geq n_k$
then
$G_n$ does not contain a cycle of length not greater than $k$. Large girth
sequences are the analogues of $L$-treeings \cite{KM}. Our first goal is to
prove the following version of Gaboriau's Theorem \cite{Gab1},
(see also Theorem 19.2  \cite{KM}).

\begin{theorem} \label{t1}
If $\cG=\{G_n\}^\infty_{n=1}$ is a large girth sequence, then
$e(\cG)=c(\cG)$.
\end{theorem}
\subsection{$\beta$-invariants}
In the proof of Theorem \ref{t1} we shall use the
$\beta$-invariants which are the analogues of the first $L^2$-Betti numbers
of measurable equivalence relations \cite{Gab2}.
First recall the notion of cycle spaces.

\noindent
Let $G(V,E)$ be a finite, simple, connected graph and $K$ be a commutative
field. Let $\varepsilon_K(G)$ be the vector space over $K$ spanned by
the edges and let $C_k(G)\subseteq \varepsilon_K(G)$, the {\bf cycle space}
be the subspace generated by the cycles of $G$. Then $\dim_K C_K(G)=
|E|-|V|+1$. Now let $\cGG$ be a graph sequence. Let $C^q_K(G_n)$ be the
space spanned by the cycles of $G_n$ of length not greater than $q$.
Here we use the usual convention that $(x,y)=-(y,x)$ and we associate
to the cycle $(x_1,x_2,\dots,x_{n},x_1)$ the vector 
$(\sum^{n-1}_{i=1}(x_i,x_{i+1})+(x_n,x_1))$.

\noindent
Set $s^q_K(\cG):=\liminf_{n\to\infty} \frac{|E(G_n)|-\dim_K C^q_K(G_n)}
{|V(G_n)|}-1\,.$ The $\beta_K$-invariant of $\cG$ is defined as
$$\beta_K(\cG):=\inf_q s^q_K(\cG)\,.$$
In Section \ref{sect2} we shall prove that if $\cG\simeq \cH$, then
$\beta_K(\cG)=\beta_K(\cH)$. This immediately shows that
$$\beta_K(\cG)+1\leq c(\cG)\,.$$
\subsection{Residually finite groups}
Let $\Gamma$ be a finitely generated group and
$$\Gamma\rhd \Gamma_1\rhd \Gamma_2\rhd\dots,\quad 
\cap^\infty_{n=1} \Gamma_n=\{1\}\,$$
be a nested sequence of finite index normal subgroups. Following
Lackenby \cite{Lac1} we define the rank gradient of the system
$\{\Gamma,\{\Gamma_n\}^\infty_{n=1}\}$
$$\mbox{rk grad}\,\{\Gamma,\{\Gamma_n\}^\infty_{n=1}\}=\limi\frac{d(\Gamma_n)}
{|\Gamma:\Gamma_n|}\,,$$
where $d(\Gamma_n)$ is the minimal number of generators for $\Gamma_n$.
In another paper \cite{Lac2}, Lackenby investigated the behaviour of the
sequence $\{\frac {d_p(\Gamma_n)}{|\Gamma:\Gamma_n|}\}^\infty_{n=1}$, where
$d_p(\Gamma_n)=\dim_{\Fp} H_1(\Gamma_n,\Fp)\,.$ 
Here we denote by $\Fp$ the finite
field of $p$ elements.
Note that $d_p(\Gamma_n)\leq d(\Gamma_n)$. The mod-$p$-homology gradient of
the system $\{\Gamma,\{\Gamma_n\}^\infty_{n=1}\}$ is defined as
$$\mbox{mod $p$-homology grad}
\,\{\Gamma,\{\Gamma_n\}^\infty_{n=1}\}=\liminf_{i\to\infty}
\frac{d_p(\Gamma_n)}
{|\Gamma:\Gamma_n|}\,.$$
Let $S$ be a symmetric generating system for $\Gamma$ and let $\cGG$ be
the graph sequence of the 
 Cayley-graphs of $\Gamma/\Gamma_n$ with respect to $S$. We have the
following theorem:
\begin{theorem} \label{t2}
$$c(\cG)-1\leq \mbox{rk grad}\,\{\Gamma,\{\Gamma_n\}^\infty_{n=1}\}\,.$$
If $\Gamma$ is even finitely presented, then
 we have the following inequality,
$$\beta_\bQ(\cG)=\beta^1_{(2)}(\Gamma)\leq \mbox{mod $p$-homology grad}
\,\{\Gamma,\{\Gamma_n\}^\infty_{n=1}\} =\beta_{\Fp}(\cG)\leq c(\cG)-1\,,$$
where $\beta^1_{(2)}(\Gamma)$ is the first $L^2$-Betti number of $\Gamma$
\cite{Lueck}. \end{theorem}
\subsection{Hyperfinite graph sequences}

One of the key notion in the theory of measurable equivalence relations
is hyperfiniteness. We introduce a similar notion for graph sequences.
We shall prove the following analogues of Proposition 22.1 and Lemma 23.2 of
\cite{KM}.
\begin{proposition}\label{hyper}
\begin{enumerate}
\item If $\cHH$ is a hyperfinite graph sequence then $c(\cH)=1$.
\item For any graph sequence $\cGG$ there exists a hyperfinite
graph sequence $\cHH$ such that $\cH\prec \cG$. 
\end{enumerate}
\end{proposition}
Finally we prove the analogue of the theorem of Connes, Feldman
and Weiss (Theorem 10.1 \cite{KM}).
\begin{theorem}\label{t3}
Let $\Gamma$ be a finitely generated residually finite group with a
nested sequence of finite index normal subgroups $\Gamma_n$, 
$\cap^\infty_{n=1} \Gamma_n=\{1\}$. Then the associated graph sequence $\cG$
is hyperfinite if and only if $\Gamma$ is amenable.
\end{theorem}
\section{$\beta$-invariants} \label{sect2}
\begin{proposition}\label{p11}
Let $\cG\simeq \cH$ be equivalent graph sequences and $K$ be a field. Then
$\beta_K(\cG)=\beta_K(\cH)$.
\end{proposition}
\begin{proof}
Suppose that $\cH\subseteq \cG$, that is for any $n\geq 1$, $E(H_n)=E(G_n)$.
Let $L>0$ be an integer such that $d_{G_n}(x,y)\leq L d_{H_n}(x,y)$.
We define a $K$-linear transformation between quotient spaces:
$$\widetilde{\phi}:\frac{\varepsilon_K(H_n)}{C^q_K(H_n)}\to 
\frac {\varepsilon_K(G_n)}{C^q_K(G_n)}$$
by extending the inclusion $\phi:E(H_n)\to E(G_n)$.
\begin{lemma} \label{l11}
$\widetilde{\phi}$ is surjective if $q>L$.
\end{lemma}
\begin{proof} Let $e=(x,y)\in E(G_n)$, then
there exists a path $P$ between $x$ and $y$,
 in $H_n$ of length not greater than $L$.
The cycle $c=P\cup e$ represents an element in $C^q_K(G_n)$ and 
$$[e]\in [c]\oplus[\widetilde{\phi}(\varepsilon_K(H_n))]\,.$$
Hence the lemma follows. \qed \end{proof}
By the lemma it follows that $s^q_K(H_n)\geq s^q_K(G_n)$ if $q>L$,
thus $\beta_K(\cH)\geq \beta_K(\cG)$.

\noindent
Now we define another $K$-linear transformation:
$$\widetilde{\psi}:\frac{\varepsilon_K(G_n)}{C^q_K(G_n)}\to 
\frac {\varepsilon_K(H_n)}{C^{qL}_K(H_n)},$$
by mapping the basis vector $e=(x,y)\in E(G_n)$ to a path in $H_n$ of length
not greater than $L$ connecting $x$ and $y$. If $e\in H_n$, then 
let $\widetilde{\psi}(e)=e$. Obviously, $\widetilde{\psi}$ is surjective
therefore $s^q_K(G_n)\geq s^{qL}_K(H_n)$ and consequently
$\beta_K(\cG)\geq \beta_K(\cH)$.

\noindent
Hence if $\cG\simeq \cH$, $\cH\subseteq \cG$ then
$\beta_K(\cG)=\beta_K(\cH)$. Now we consider the general case, where
$\cH$ is an arbitrary graph sequence such that $\cH\simeq \cG$.
Then let $\cJ=\cG\cup \cH$, that is
$V(J_n)=V(G_n), E(J_n)=E(G_n)\cup E(H_n)$.\
Clearly, $\cJ\simeq \cG\simeq \cH$ and $\cH\subseteq \cJ$, $\cG\subseteq \cJ$.
Thus by our argument above, $\beta_K(\cH)=\beta_K(\cJ)=\beta_K(\cG)$\qed
\end{proof}
\begin{proposition} \label{p14}
Let $\cGG$ be a graph sequence then 
$$\beta_\bQ(\cG)\leq \beta_{\Fp}(\cG)\leq c(\cG)-1\,.$$
\end{proposition}
\begin{proof}
Let $\cH\simeq \cG$, then $\beta_K(\cG)=\beta_K(\cH)\leq e(\cH)-1\,$. 
Therefore $\beta_K(\cG)\leq c(\cG)-1\,.$
\begin{lemma}\label{l14}
$\dim_{\bQ} C^q_{\bQ}(G_n)\leq \dim_{\Fp} C^q_{\Fp}(G_n)\,.$
\end{lemma}
\begin{proof}
Let $c^q_n$ be the number of cycles in $G_n$ of length not greater than $q$.
Let $\rho_{\bZ}:\bZ^{c^q_n}\to \bZ^{|E(G_n)|}$ be the homomorhism that maps
$\oplus^{c^q_n}_{i=1} s_i$ to $\sum^{c^q_n}_{i=1} s_i [c_i]$, where
$s_i\in\bZ$ and $[c_i]$ is the integer vector generated by the $i$-th cycle
$c_i$. Similarly, we define
$\rho_{\Fp}:\Fp^{c^q_n}\to \Fp^{|E(G_n)|}$. 
Let $\pi_1:\bZ^{c^q_n}\to\Fp^{c^q_n}$, $\pi_2:
\bZ^{|E(G_n)|}\to\Fp^{|E(G_n)|}$
be the residue class maps. Then $\pi_2\circ \rho_{\bZ}=\rho_{\Fp}\circ
\pi_1\,.$
Therefore,
$$\rank_{\bZ}\im \rho_{\bZ}\geq \dim_{\Fp} \im \rho_{\Fp}\,.$$ 
Clearly, $\rank_{\bZ}\im \rho_{\bZ}=\dim_{\bQ} C^q_{\bQ} (G_n)\,$
and 
$\dim_{\Fp}\im \rho_{\Fp}=\dim_{\Fp} C^q_{\Fp} (G_n)\,.$
Thus our lemma follows. \qed \end{proof}
By our lemma, $\beta_{\bQ}(\cG)\leq \beta_{\Fp}(\cG)$ hence
we finish the proof of our proposition. \qed \end{proof}
\begin{question}
Does there exist a graph sequence $\cG$ for which
$\beta_{\bQ}(\cG)\neq \beta_{\Fp}(\cG)$ or
$\beta_{\Fp}(\cG)\neq c(\cG)-1$ ?
\end{question}
Finally we prove Theorem \ref{t1}.

\noindent
\begin{proof}
Let $\cGG$ be a large girth graph sequence. Then by definition
$\beta_K(\cG)=e(\cG)-1\,.$ That is $e(\cG)-1\leq c(\cG)-1$, hence our
theorem follows. \qed \end{proof}
\section{Residually finite groups}
The goal of this section is to prove Theorem \ref{t2}. 
Let $\Gamma$ be a finitely generated residually finite group with
a not necessarily symmetric generating system $S$. Let
$\Gamma\rhd \Gamma_1\rhd \Gamma_2\rhd\dots,\quad 
\cap^\infty_{n=1} \Gamma_n=\{1\}\,$
be a nested sequence of finite index normal subgroups and
$\cGG$ be a graph sequence, where $G_n$ is the (left) Cayley-graph
of the finite group $\Gamma/\Gamma_n$ with respect to $S$. Note that if
$S'$ is another generating system and $\cHH$ is the associated graph
sequence then $\cH\simeq \cG$.
\begin{proposition} \label{p17}
$c(\cG)-1\leq \mbox{rk grad\,} \{\Gamma, \{\Gamma_n\}^\infty_{n=1}\}\,.$
\end{proposition}
\begin{proof}
First note that by the Reidemeister-Schreier theorem the groups $\Gamma_n$
are finitely generated as well \cite{MKS}, moreover if $T$ is
finite generating system of $\Gamma_n$, then
$$d_{G^{\Gamma_n}_T}(x,y)\leq L d_{G^\Gamma_S}(x,y)$$ for any 
$x,y\in \Gamma_n$,
where $G^\Gamma_S$ resp. $G^{\Gamma_n}_T$ are the
Cayley-graphs with respect to $S$ resp. to $T$, and the Lipschitz constant
$L$ depends only on $S$ and $T$.
\begin{lemma}
\label{l18}
For any $k\geq 1$, 
$$\frac{d(\Gamma_k)}{|\Gamma:\Gamma_k|}+1\geq c(\Gamma)\,.$$
\end{lemma}
\begin{proof}
We use an idea resembling an argument in the proof of Theorem 21.1 \cite{KM}.
Let $T$ be a generating system of $\Gamma_k$ of minimal number of
generators. For simplicity we suppose that $T\subset S$.
Consider the following graph sequence $\cHH$, $V(H_n)=\Gamma/\Gamma_n$.
If $n\leq k$, let $H_n=G_n$. Set $S_n=\Gamma_k/\Gamma_n$ and let
$H'_n$ be the Cayley-graph of $S_n$ with respect to $T$.
Now enumerate the vertices of
$V(H_n)\backslash S_n$, $\{x_1, x_2\dots, x_{r_n}\}$. For each $x_i$
consider the set of shortest paths in $G_n$ from $x_i$ to the set $S_n$. Pick
the minimal path with respect to the lexicographic ordering. The edges
of $H_n$ shall consist of $H'_n$ and the edges of the minimal paths.
Define a map $\pi:V(H_n)\to S_n$ the following way.
For each $x_i\in V(H_n)\backslash S_n$ let $\pi(x_i)\in S_n$ be the endpoint
of the minimal path from $x_i$ to $S_n$ and let $\pi(x)=x$ if $x\in S_n$.
By the lexicographic minimality, the union of the paths form a subforest
in $G_n$ having exactly $|V(H_n)\backslash S_n|$ edges.

\noindent
We claim that $\cH\simeq \cG$. Since $\cH\subset \cG$, we only need to
prove that $\cG\prec\cH$. Let $n>k$, $x,y\in V(G_n)$.
Consider the shortest $G_n$-path from $x$ to $y$, 
$\{x_0,x_1,\dots x_l\}$, $x_0=x$, $x_l=y$. Let us consider the
sequence of vertices $\{\pi(x_0), \pi(x_1),\dots \pi(x_l)\}$.

\noindent
Let $y_1,y_2,\dots, y_{|\Gamma:\Gamma_k|}$ be a set of coset-representatives
with respect to $\Gamma_k$. Let $t$ be the maximal word-length of the
representatives with respect to $S$. Then  $d_{G_n}(\pi(x),x)\leq t$
holds for any $x\in V(G_n)$. Therefore,
$d_{G_n}(\pi(x_i), \pi(x_{i+1}))\leq 2t+1\,.$
That is $d_{H_n}(\pi(x_i),\pi(x_{i+1})\leq L(2t+1)$, where
$L$ is the Lipschitz-constant defined before stating our lemma.
Consequently,
$$d_{H_n}(x,y)\leq L(2t+1) d_{G_n}(x,y)$$ 
and therefore $\cH\simeq\cG$. 

\noindent
For the edge number of $\cH$ we have
$$e(\cH)=\liminf_{n\to\infty}\frac
{|\Gamma:\Gamma_n|-|\Gamma_k:\Gamma_n| +|E(H_n'|}
{|\Gamma:\Gamma_n|}\,.$$
The vertex degrees of $H'_n$ are not greater than $2|T|=2d(\Gamma_k)$, also
$|S_n|=|\Gamma_k:\Gamma_k|\,.$ Thus
$$c(\cG)\leq e(\cH)\leq \frac{d(\Gamma_k)}{|\Gamma:\Gamma_k|}+1\,.$$
Hence the lemma follows. \qed 
\end{proof}
Proposition \ref{p17} is a straightforward consequence of Lemma \ref{l18}
\qed
\end{proof}
Let $\{\Gamma,\{\Gamma_n\}^\infty_{n=1}\},S,\cG$ be as above. Moreover suppose
that $\Gamma$ is finitely presented. This means that if $\Theta
:\cF_S\to
\Gamma$ is the natural map from the free group generated by $S$ to $\Gamma$
then
$\ke \cF_S$ is generated by the relations $\{R_1,R_2,\dots,R_l\}$ as normal
subgroup that is if $\Theta(\underline{w})=1$ then
$$ \underline{w}=\prod^{r_{\underline{w}}}_{j=1} 
\gamma_{j}R_{i_j}\gamma_j^{-1}\,,\quad
\gamma_j\in\cF_s\,.$$
Let $\widetilde{\Sigma}$ be the usual covering $CW$-complex constructed
from $\{R_i\}^l_{i=1}$, the $1$-skeleton of $\widetilde{\Sigma}$ is
the Cayley-graph of $\Gamma$ and for 
each $\gamma\in\Gamma$ and $1\leq i \leq l$, we add a $2$-cell
$\sigma_{\gamma,i}$ such that
$$\partial\sigma_{\gamma,i}=\sum^{s_i}_{j=1}
(\underline{w}_j\gamma,\underline{w}_{j-1}\gamma)\,,$$
where $R_i=a_{s_i}a_{s_i-1}\dots a_2a_1\,,$ $\underline{w}_j=
a_ja_{j-1}\dots a_2 a_1$, $\underline{w}_0=1$.
Then $\widetilde{\Sigma}$ is simply connected with a natural $\Gamma$-action.
Clearly, $\pi_1(\widetilde{\Sigma}/\Gamma_n)=\Gamma_n$. 
Recall that the group homology space $H_1(\Gamma_n,K)$ is
isomorphic to the $CW$-homology space $H_1(\widetilde{\Sigma}/\Gamma_n,K)$.
\begin{lemma}
\label{l24}
$\limn\frac{\dim_K H_1(\widetilde{\Sigma}/\Gamma_n,K)}{|\Gamma:\Gamma_n|}=
\beta_K(\cG)$.
\end{lemma}
\begin{proof}
Consider the homology complex
$$C_2(\widetilde{\Sigma}/\Gamma_n,K)\stackrel{\partial_2}{\to}
C_1(\widetilde{\Sigma}/\Gamma_n,K)\stackrel{\partial_1}{\to}
C_0(\widetilde{\Sigma}/\Gamma_n,K)\,.$$
Observe that 
$$C_1(\widetilde{\Sigma}/\Gamma_n,K)\simeq\varepsilon_K(G_n)\,\quad
\dim_K C_0(\widetilde{\Sigma}/\Gamma_n,K)=|V(G_n)|\,.$$
Let $r$ be the maximal word-length of a relation $R_i$. Then
$\partial_2(C_2(\widetilde{\Sigma}/\Gamma_n,K))$ is generated by cycles
of length at most $r$. On the other hand for any $q>r$, the $q$-cycles
are in $\partial_2(C_2(\widetilde{\Sigma}/\Gamma_n,K))$ if $n$ is large
enough. 

\noindent
Therefore $C^q_K(G_n)=\partial_2(C_2(\widetilde{\Sigma}/\Gamma_n,K))$
if $n$ is large enough. Consequently
$$s^q_K(\cG)=\liminf_{n\to\infty}
\frac{|E(G_n)|-\dim_K \partial_2(C_2(\widetilde{\Sigma}/\Gamma_n,K))-|V(G_n)|}
{|\Gamma:\Gamma_n|}\,.$$
On the other hand,
$$\frac{\dim_K H_1(\widetilde{\Sigma}/\Gamma_n,K)}{|\Gamma:\Gamma_n|}=
\frac{\dim_K \ke \partial_1 -\dim_K \im \partial_2}{|\Gamma:\Gamma_n|}=$$
$$=
\frac{|E(G_n)|-\dim_K \partial_2(C_2(\widetilde{\Sigma}/\Gamma_n,K))-
|V(G_n)|+1}{|\Gamma:\Gamma_n|}\,.$$
Hence the lemma follows. \qed \end{proof}

Now we prove the second part of Theorem \ref{t2}. 
\begin{proposition}
\label{p26}
Let $\Gamma$ be a finitely presented residually finite group,
$\{\Gamma,\{\Gamma_n\}^\infty_{n=1}\},S,\cG$ be as above. Then
$$\beta_\bQ(\cG)=\beta^1_{(2)}(\Gamma)\leq \mbox{mod $p$-homology grad}
\,\{\Gamma,\{\Gamma_n\}^\infty_{n=1}\} =\beta_{\Fp}(\cG)\leq c(\cG)-1\,,$$
where $\beta^1_{(2)}(\Gamma)$ is the first $L^2$-Betti number of $\Gamma$
\cite{Lueck}.
\end{proposition}
\begin{proof}
By Lemma \ref{l24}
$\beta_{\Fp}(\cG)=\mbox{mod $p$-homology grad}
\,\{\Gamma,\{\Gamma_n\}^\infty_{n=1}\}$. Also,
$$\beta_\bQ(\cG)=\liminf_{n\to\infty} \frac{\dim_{\bQ}
H_1(\widetilde{\Sigma}/\Gamma_n,\bQ) }{|\Gamma:\Gamma_n|}\,.$$
By the Approximation Theorem of L\"uck
$$\lim_{n\to\infty} \frac{\dim_{\bQ}
H_1(\widetilde{\Sigma}/\Gamma_n,\bQ) }{|\Gamma:\Gamma_n|}=
\beta^1_{(2)}(\Gamma)\,.$$
Hence our proposition follows. \qed \end{proof}
\begin{question}
\begin{enumerate}
\item Does there exist a finitely presented residually finite group $\Gamma$
and a system $\{\Gamma,\{\Gamma_n\}^\infty_{n=1}\}$ such that
$$\beta^1_{(2)}(\Gamma)\neq \mbox{mod $p$-homology grad}
\,\{\Gamma,\{\Gamma_n\}^\infty_{n=1}\}\quad\mbox{or} $$ $$
\mbox{mod $p$-homology grad}
\,\{\Gamma,\{\Gamma_n\}^\infty_{n=1}\}\neq c(\cG)-1\, ?$$
\item Does there exist a finitely generated residually finite group $\Gamma$
and a system $\{\Gamma,\{\Gamma_n\}^\infty_{n=1}\}$ such that
$$c(\cG)-1\neq \mbox{rk grad}\,\{\Gamma,\{\Gamma_n\}^\infty_{n=1}\}\, ?$$
\end{enumerate}
\end{question}
\section{Hyperfinite graph sequences}
We say that a graph sequence $\cGG$ is {\bf hyperfinite}
if for any $\e>0$ there exists $K_\e>0$ and a sequence of partitions of
the vertex sets $V(G_n)$
$$A^n_1\cup A^n_2\cup\dots\cup A^n_{k_n}=V(G_n)\,$$
such that
\begin{itemize}
\item For any $n\geq 1$, $1\leq i \leq k_n$, $|A^n_i|\leq K_\e$.
\item If $E^\e_n$ is the set of edges $(x,y)\in E(G_n)$ such that
$x\in A_i$, $y\in A_j$, $x\neq y$, then
$$\liminf_{n\to\infty} \frac{|E^\e_n|}{|V(G_n)|}\leq \e\,.$$
\end{itemize}
Now we prove Proposition \ref{hyper}.

\noindent
\begin{proof}
Suppose that $\cGG$ is hyperfinite. Let $\cH^\e=\{H_n^\e\}^\infty_{n=1}$ be
 the following
graph sequence. The vertex set of $H^\e_n$ is $V(G_n)$, $E(H^\e_n)$ is
the union of $E^\e_n$ and a spanning tree for each connected component
of the graphs spanned by the vertices of $A^n_i$, $1\leq i \leq k_n$.
Clearly, $\cH^\e\simeq \cG$ and $|E(H^\e_n)|\leq |E^\e_n|+|V(G_n)|$ thus
$e(\cH^\e)\leq 1+\e$. Therefore $c(\cG)=1\,.$

\noindent
Now we show that for any graph sequence $\cGG$, $\cHH$ is hyperfinite
where $H_n$ is a spanning tree of $G_n$. We actually show that a sequence
of trees $\cT=\{T_n\}^\infty_{n=1}$ is always hyperfinite.
Let $q$ be an integer and consider a maximal $q$-net $L^q_n\subset V(T_n)$.
That is if $x\neq y\in L^q_n$ then $d_{T_n}(x,y)\geq q$ and
for any $z\in V(T_n)$ there exists $x\in L^q_n$ such that $d_{T_n}(x,z)\leq
q$.
Now for each $x\in V(T_n)$ let $\pi(x)$ be one of the vertices $y\in L^q_n$
closest to $x$. Then $\cup_{y\in L^q_n} \pi^{-1}(y)$ is a partition
of $V(T_n)$.
Clearly $|\pi^{-1}(y)|\geq q$ for any $y\in L^q_n$. Obviously the $T_n^y$
subgraph spanned by the vertices in $\pi^{-1}(y)$ is connected. Thus
$$|E^\e_n|\leq |V(T_n)|-(|V(T_n)|-|L^q_n|)\,.$$
Here we used the fact that a connected graph has at least as many edges
as the number of its vertices minus one.
Obviously, $|L^q_n|\leq\frac{|V(T_n)|}{q}$, therefore
$$\limn \frac{|E^\e_n|}{|V(T_n)|}\leq \frac{1}{q}\,.$$
Consequently, the graph sequence $\cT$ is indeed hyperfinite. \qed
\end{proof}
Finally, we prove Theorem \ref{t3}

\noindent
\begin{proof}
First let $\Gamma$ be a residually finite non-amenable group with 
a symmetric generating system $S$ and a nested sequence of finite
index normal subgroups 
$\Gamma\rhd \Gamma_1\rhd \Gamma_2\rhd\dots,\quad 
\cap^\infty_{n=1} \Gamma_n=\{1\}\,$. Let $G_n$ be the Cayley-graph
of $\Gamma/\Gamma_n$ with respect to $S$ and $G^\Gamma_S$ be the Cayley-graph
of the group $\Gamma$. Since $\Gamma$ is non-amenable, it has no
\Fo-exhaustion, consequently there exists a real number $\delta>0$
such that for each finite subset $F\subset \Gamma$ the number of edges from 
$F$ to the complement of $F$ is at least $\delta|F|$.
Fix an integer $m>0$. If $n$ is large enough then for any subset 
$M\subset \Gamma/\Gamma_n$, $|M|\leq m$ the number of edges from $M$
to its complement must be at least $\delta M$. This follows easily form the
fact that for any $r\geq 0$, the $r$-balls in $G_n$ and in $G^\Gamma_S$
are isometric. This implies that $\cG$ is not hyperfinite.

\noindent
Now let $\Gamma, \{\Gamma_n\}^\infty_{n=1}, S,\cG$ be as above, but 
 let $\Gamma$ be amenable. The following lemma is a straight consequence
of Theorem 2. \cite{Elek}
\begin{lemma} \label{l33}
For any $\om>0$, there exist $L_\om>0$, $M_\om>0$ and
a sequence of family of subsets
$$\{W^i_n\}^{k_n}_{i=1},\quad W^i_n\subset V(G_n)\quad\mbox{if}\quad
n\geq M_\om$$
such that for any $1\leq i \leq k_n$,
\begin{itemize}
\item $|W^i_n|\leq L_\om\,.$
\item $|W^i_n\backslash \cup^{k_n}_{j\neq i} W^j_n|\geq (1-\omega)|W^i_n|\,.$
\item The number of edges from $W^i_n$ to its complement is at most
$\om|W^i_n|$.
\end{itemize}
and
\begin{itemize}
\item $|\cup^{k_n}_{i=1} W^i_n|\geq (1-\om)|V(G_n)|$.
\end{itemize}
\end{lemma}
Now let $Z^i_n=W^i_n\backslash \cup^{k_n}_{j\neq i} W^j_n\,$ and consider
the partition of $V(G_n)$
$$V(G_n)=\bigcup^{k_n}_{i=1} Z^i_n \cup \bigcup^{l_n}_{j=1} T^i_n\,,$$
where $T^i_n$ are arbitrary subsets of size at most $L_\om$.
Let $E^\om_n$ be the set of edges $(x,y)\in G_n$ such that their
endpoints belong to different subsets in the partition. There are
three kinds of edges in $E^\om_n$:
\begin{itemize}
\item Edges with an endpoint in $T^i_n$. The number of such
edges is at most $2|S| (1-(1-\om)^2)|V(G_n)|\,.$
\item Edges from $Z^i_n$ to the complement of $W^i_n$, for some $1\leq i \leq
  k_n\,.$ The number of such edges is at most $2|S| \om(1-\om)^{-1}
  |V(G_n)|\,.$
\item Edges from $Z^i_n$ to $W^i_n\backslash Z^i_n$ for some $1\leq i \leq
  k_n$.
The number of such edges is at most  $2|S| \om(1-\om)^{-1}
  |V(G_n)|\,.$
\end{itemize}
Hence
$$\liminf_{n\to \infty} \frac{|E^\om_n|}{|V(G_n)|}\leq 
2|S|((1-(1-\om)^2)+2  \om(1-\om)^{-1})\,.$$
Therefore $\cG$ is hyperfinite. \qed
\end{proof}

\begin{address}
G\'abor Elek \\
The Alfred Renyi Institute of the Hungarian Academy of Sciences
H-1364, P.O.B 127
Hungary
\email{elek@renyi.hu}
\end{address}

\end{document}